\documentclass[letterpaper,conference,10pt]{ieeeconf}
\IEEEoverridecommandlockouts
\usepackage{stmaryrd} 
\SetSymbolFont{stmry}{bold}{U}{stmry}{m}{n} 
\usepackage[dvipsnames]{xcolor}

\usepackage[sort,noadjust]{cite}
\usepackage{amsmath,amssymb,amsfonts}

\usepackage{nicefrac}

\usepackage{soul}
\sethlcolor{cyan}
\usepackage{textcomp}
\def\BibTeX{{\rm B\kern-.05em{\sc i\kern-.025em b}\kern-.08em
    T\kern-.1667em\lower.7ex\hbox{E}\kern-.125emX}}
\usepackage{dsfont}
 
\usepackage[T1]{fontenc}

\usepackage{graphicx}

\usepackage{amssymb,amsmath}

\usepackage{xparse}
\usepackage{amsmath}
\usepackage{bm}
\usepackage{mathtools}

\usepackage{tabularx}

\usepackage{cases}
\usepackage{diagbox}
\usepackage{algorithm}

\usepackage{algpseudocode}

\usepackage[hidelinks,
      breaklinks=true,
]{hyperref}

\usepackage{booktabs}
\usepackage{multirow}

\usepackage{tabularx}
\usepackage{tcolorbox}

\newcommand\numberthis{\addtocounter{equation}{1}\tag{\theequation}}

\newcommand{\Rl}[2]{\ensuremath{\mathbb{R}^{#1}_{#2}}}   
\newcommand{\R}{\ensuremath{\mathbb{R}}}

\newcommand{\Rlo}{\ensuremath{\mathbb{R}_{\geq 0}}}

\newcommand{\Zo}{\ensuremath{\mathbb{Z}_{\geq 0}}}
\newcommand{\Zp}{\ensuremath{\mathbb{Z}_{> 0}}}

\definecolor{bleucit}{rgb}{0.2,0.4,0.6} 

\definecolor{blue_cv}{rgb}{0.09,0.35,0.78}

\newcommand{\Argmin}{\ensuremath{\text{argmin}\,}}
\DeclareMathOperator*{\argmin}{\Argmin}

\DeclareMathOperator{\Tr}{Tr}

\newcommand{\dom}{\ensuremath{\text{dom}\,}}

\newtheorem{prop}{Proposition}

\newtheorem{lem}{Lemma}

\newtheorem{thm}{Theorem}

\newtheorem{rem}{Remark}
\newtheorem{sass}{Standing Assumption}

\usepackage[all]{hypcap}

\usepackage{svg}   
\usepackage{subcaption}

\title{\LARGE \bf
An optimistic planning algorithm for switched discrete-time LQR}

\author{\vspace{-1ex}\thanks{Work funded by the ANR under grant OLYMPIA ANR-23-CE48-0006, by the ARC under the Discovery Project DP210102600 and DP250100300\textcolor{black}{, and by the CINEA under grant SeaClear2.0,  No 101093822.}} Mathieu Granzotto\thanks{Mathieu Granzotto and Dragan Nešić are with the Department of Electrical and Electronic Engineering, University of Melbourne, Parkville, VIC 3010, Australia (e-mail: \{mgranzotto, dnesic\}@unimelb.edu.au).}, Romain Postoyan\thanks{Romain Postoyan and Jamal Daafouz  are with the Universit\'e de Lorraine, CNRS, CRAN, F-54000 Nancy, France (e-mails: \{name.surname\}@univ-lorraine.fr). J. Daafouz is also with IUF.}, Dragan Nešić,  Jamal Daafouz and Lucian Buşoniu\thanks{Lucian Buşoniu is with the  Technical University of Cluj-Napoca, 400114 Cluj-Napoca, Romania (e-mail: lucian.busoniu@aut.utcluj.ro) \textcolor{black}{and is a  Corresponding
Member of the Romanian Academy.}}} 

\usepackage{xspace}
\newcommand{\OPlin}{\textsc{TROOP}\xspace}

\begin{document}

\bstctlcite{BSTcontrol}
 
\maketitle
\thispagestyle{empty}
\pagestyle{empty}
\setlength{\textfloatsep}{0pt}
\allowdisplaybreaks


\newcommand{\assname}{Assumption}
\newcommand{\corname}{Corollary}
\newcommand{\thmname}{Theorem}
\newcommand{\propname}{Proposition}
\newcommand{\lemname}{Lemma}
\newcommand{\algorithmname}{Algorithm}
\newcommand{\remname}{Remark}
\newcommand{\defnname}{Definition}
\renewcommand{\sectionautorefname}{Section}
\renewcommand{\subsectionautorefname}{Section}
\renewcommand{\subsubsectionautorefname}{Section}
\makeatletter
\newcommand{\sassname}{SA\@gobble}
\makeatother
\makeatletter
\newcommand{\tcname}{C\@gobble}
\makeatother



\begin{abstract}
We introduce \OPlin, a tree-based Riccati optimistic online planner, that is designed to generate   near-optimal control laws   for discrete-time switched linear systems with switched quadratic costs. The key challenge that we address is balancing computational resources against control performance, which is important as constructing near-optimal inputs often requires substantial amount of computations. \OPlin addresses this trade-off by adopting an online best-first search strategy inspired by $A^\star$, allowing for efficient estimates of the optimal value function.
The control laws obtained  guarantee both near-optimality and stability properties for the closed-loop system. These properties depend on the planning depth, which determines how far into the future the algorithm  explores and is closely related to the amount of computations. \OPlin thus strikes a balance between computational efficiency and control performance, which is illustrated by numerical simulations on an example.
\end{abstract}

\section{Introduction}
 
Consider the problem of regulating a vehicle's  speed while minimizing fuel consumption. Achieving this requires selecting the best gear (a discrete input) and throttle level (a continuous input), leading to an optimal hybrid-valued control scheme, specifically a mixed-integer problem. In view of the combinatorial nature of mixed-integer problems, it is   key to balance the computational resources available against desired closed-loop system properties. Specifically, we typically want to obtain a cost close to the minimum one (near-optimality) while ensuring that the selected inputs drive the system to a desired operating region (stability). Even when the cost function is quadratic and the system dynamics are  linear for each value of the discrete input, this remains computationally very challenging \cite{antunes_linear_2017}. This latter setting describes the switched LQR problem: to determine a sequence of discrete linear modes as well as a sequence  of continuous input vectors which, together, minimize switched quadratic costs along the solution of the closed-loop system. Switched  problems arise in a wide range of important applications, including power systems and power electronics (like switching supplies), automotive control (like gear selection above), aircraft and air traffic control, computer disks, and
network and congestion control   \cite{ge_switched_2005}. In this work, we study the discrete-time switched LQR problem. 

 The switched LQR problem has received much attention, see, e.g., \cite{weizhang_Value_2009,zhang_InfiniteHorizon_2012,hou_Improved_2023,zhao_Optimal_2020}. The bulk of existing approaches rely on approximating the optimal   value function for any state by a piecewise quadratic function, given by sets of positive semi-definite  matrices \cite{antunes_linear_2017,hou_Improved_2023,zhang_InfiniteHorizon_2012,beumerComplexityBoundedRelaxedDynamic2023a}. These matrices are then used to construct  the discrete and   continuous inputs. However, when generating   tighter  approximations of  the optimal value function, these sets of matrices generally become  exponentially large in the number of discrete modes denoted $M$  \cite{weizhang_Value_2009}. There is a need to provide algorithms for the  switched LQR problem that trade control performance with  computational complexity. Indeed, while several   computational strategies are available, such as   employing pruning strategies \cite{zhang_InfiniteHorizon_2012,hou_Improved_2023}, exploiting a ``good'' initial choice of cost matrices \cite{antunes_linear_2017}, or establishing a unique periodic sequence by nesting Linear Matrix Inequalities (LMIs) \cite{ji-woongleeInfiniteHorizonJointLQG2009}, these approaches often fail to balance near-optimality with computational efficiency.    Most relevantly, the approach of \cite{weizhang_Value_2009} characterizes a ``breadth-first'' exploration of all available modes of the discrete inputs, expanding all sequences of a given length. These solve a finite-horizon switched LQR problem with   ever-increasing horizons  (by increasing the depth of the exploration tree), thus producing ever-improving approximations of the optimal cost.  Although calculated offline, this computationally intensive approach can still lead to controllers requiring complex online searches in large lists, with potentially unsatisfactory closed-loop
performance, even when employing aggressive pruning strategies, see \cite[Section V]{antunes_linear_2017}.

In contrast, we propose an online planning algorithm that performs a “best-first” exploration of modes, prioritizing sequences with lower costs based on the current state. We call this new (near-)optimal control algorithm \OPlin (tree-based Riccati optimistic online planner).  Drawing inspiration from efficient planning algorithms like $A^\star$ and optimistic planning \cite{Hren-Munos-2008optimistic,Busoniu-Munos-12,granzotto_stable_2022}, the proposed approach typically eliminates many branches from exploration. That is, instead of requiring $\mathcal{O}(M^d)$ iterations  for a tree of horizon (or depth) $d$, where $M$ is the number of discrete inputs, we require $\mathcal{O}(c^d)$ with $c\leq M$, with $c$ potentially close to one in the best-case scenario. By ignoring high-cost input sequences, our method produces finite-horizon optimal control laws similar to those in \cite{weizhang_Value_2009} but with significantly reduced computational costs. This enables the efficient generation  of long input sequences with desirable closed-loop properties.

Planning  algorithms typically require an enumerable  number of possible future states    to generate the exploration tree, while the switched LQR problem involves inputs taking values in a continuous space. A key insight  of the switched LQR problem is that, given any fixed and finite sequence of discrete inputs, the optimal continuous input is given by a time-varying Ricatti equation, see, e.g.,  \cite{lincoln_lqr_2002,lincoln_efficient_2000}. This  removes the need to explicitly consider the continuous input or future states for planning, and instead we can concentrate on developing   a tree of discrete input sequences and their  associated    Ricatti equations. This is in contrast to related optimistic planning algorithms \cite{busoniuOptimisticPlanningContinuousaction2013,busoniuContinuousactionPlanningDiscounted2018}, that search  over continuous inputs,  and \cite{Hren-Munos-2008optimistic,Busoniu-Munos-12,granzotto_stable_2022}, which explore a tree of future states, which are here unfeasible. \textcolor{black}{In other words, by leveraging the switched linear structure, we significantly reduce the computational costs associated with handling continuous inputs.} 
 
 A feature of \OPlin is that the matrices computed   are state-independent, which implies that the core and computationally intensive calculations can be saved in a cache, to be later re-used in the calculation of costs for longer sequences or other initial states. Furthermore, we  allow for an initial cost function that lower bounds the optimal value function, which accelerates exploration and reduces online calculations. This is done by exploiting the structure of the switched LQR problem, providing an additional benefit of our method when compared to related optimistic planning algorithms, in particular \cite{granzotto_stable_2022}.

To achieve the desired properties of \OPlin, we rely on assumptions commonly used in the switched LQ literature. Specifically, we assume as in \cite{zhang_Study_,zhang_InfiniteHorizon_2012,weizhang_Value_2009,antunes_linear_2017,hou_Improved_2023}  that the system can be driven to the origin with bounded quadratic cost and that the stage-costs are positive definite with respect to the origin. We also impose a   condition on the initial cost function which can always be enforced, made to guarantee a monotonicity property of the estimates of the optimal value function.  Under these assumptions, \OPlin offers tunable near-optimality and global exponential stability, both dependent on the planning horizon. Longer horizons yield tighter value function approximations and improved performance at the cost of more computations.

In sum,  our approach ties together the rationales of \cite{zhang_InfiniteHorizon_2012} and \cite{antunes_linear_2017}, employing an online yet efficient optimistic planning algorithm based on \cite{Busoniu-Munos-12,granzotto_stable_2022}, leveraging known bounds of the optimal value function to reduce controller complexity. Additionally, the switching sequence and feedback gain are optimized based on the current state, offering tighter near-optimality compared to offline sequence planning methods like \cite{ji-woongleeInfiniteHorizonJointLQG2009}.   Overall, \OPlin offers a reliable framework for applications requiring both near-optimality and stability.

\section{Preliminaries} \label{sec:notation}
\subsection{Notation} \label{subsec:notation}
Let $\R$ be the set of real numbers, 
$\Rlo:=[0,\infty)$,  
$\Zp:=\{1,2,\ldots\}$, $\Zo:=\{0,1,2,\ldots\}$.
The notation $(x,y)$ stands for $[x^{\top},\,y^{\top}]^{\top}$, where $x\in\Rl{n}{}$,  $y\in\Rl{m}{}$ and $n,m\in\Zp$.  We denote the Euclidean norm by $|\cdot|$.  
The identity map  from $\Rlo$ to $\Rlo$ is denoted by $\mathbb{I}_{\geq0}$, and the zero map from $\Rlo$ to $\{0\}\subset\R$ by $\bm{0}$.    
For a (possibly infinite) sequence $\bm{u}=(u_0,\ldots,u_{d-1})$ of $d\in\Zo\cup\{\infty\}$  elements that take values in $\R^{n}$ with $n\in\Zp$, we define $\bm{u}_{\llbracket a, b\rrbracket}:=(u_a,\ldots,u_b)$ with $0\leq a\leq b < d$, and furthermore define   $\bm{u}_{\llbracket \cdot, b\rrbracket}:=\bm{u}_{\llbracket 0, b \rrbracket}$, $\bm{u}_{\llbracket a,\cdot \rrbracket}:=\bm{u}_{\llbracket a, d-1 \rrbracket}$ and $\bm{u}_{\llbracket d,\cdot\rrbracket}:=\bm{u}_{\llbracket \cdot,-1 \rrbracket}:=\emptyset$ by convention.  The trace of   real square matrix $A$ is denoted by $\text{tr}(A)$, while $\text{diag}$ refers to constructing a block diagonal matrix from its arguments.  The smallest and the largest eigenvalue of a real, symmetric  matrix $A$ are denoted $\underline{\lambda}(A)$ and $\overline{\lambda}(A)$, respectively. We write $A\succeq0$  when real, square matrix $A$ is symmetric and positive semi-definite, and $A\succ 0$ when it is symmetric and positive definite. We write
$A\succeq B$ when $A-B\succeq0$,  and  $A\succ B$ when $A-B\succ0$, for any real, symmetric matrices $A,B$. For any matrix $A\succeq0$, $\sqrt{A}$ denotes a positive semi-definite matrix such that $\sqrt{A}^\top\sqrt{A}=A$. A set-valued map $F: \mathcal{X}\rightrightarrows\mathcal{Y}$ is a multivalued function that maps elements of   $\mathcal{X}$ to subsets of $\mathcal{Y}$. Given a set-valued map $S:\R^n \rightrightarrows\R^m$, a selection of $S$ is a single-valued mapping $s:\dom S \to \R^m$ such that $s(x)\in S(x)$ for any $x\in \dom S$. For the sake of convenience, we write $s\in S$ to denote a selection $s$ of $S$. Function composition is denoted by $\circ$.

\subsection{Problem statement}

Consider the system 
\begin{equation}
    x^+ = A_i x + B_i u \label{eq:sys-lin-mixed}
\end{equation}
where $x\in\R^{n_x}$, $u\in\R^{n_u}$ and $i\in\mathbb{M}:=\{1,\ldots,M\}$ are  the state, the continuous control input   and  the discrete control input, respectively, with $n_x,n_u,M\in\Zp$. Given $k\in\Zo$ and a pair of sequences $\bm{u}_{\llbracket \cdot, k \rrbracket},\bm{i}_{\llbracket  \cdot, k \rrbracket}$ with elements in $\R^{n_u}$ and $\mathbb{M}$, respectively,   we denote the solution of \eqref{eq:sys-lin-mixed} initialized at time $0$ and state $x$ by $\phi$, which thus satisfies  $\phi(k,x,\bm{u}_{\llbracket  \cdot, k \rrbracket},\bm{i}_{\llbracket  \cdot, k \rrbracket}):=A_{i_{k-1}} \phi(k-1,x,\bm{u}_{\llbracket  \cdot,k-1 \rrbracket},\bm{i}_{\llbracket \cdot,k-1 \rrbracket})+B_{i_{k-1}} u_{k-1}$. We   define the   cost associated with a solution to system~\eqref{eq:sys-lin-mixed} initialized at $x$  with sequence of continuous and discrete input $\bm{u}$ and $\bm{i}$ by
\begin{equation}\label{eq:J}J (x,\bm{u},\bm{i}):=\sum_{k=0}^{\infty} \ell(\phi(k,x,\bm{u}_{\llbracket \cdot, k \rrbracket},\bm{i}_{\llbracket \cdot, k \rrbracket}),u_k,i_k)\end{equation}
where  $\ell(x,u,i):=x^\top Q_i x + u^\top R_i u$ for any $x\in\R^{n_x}$, $u\in\R^{n_u}$ and $i\in\mathbb{M}$.  We assume that matrices $A_i,B_i,Q_i,R_i$ are real and of conformable dimensions, that $Q_i,R_i\succ0$ for all $i\in\mathbb{M}$, and that system \eqref{eq:sys-lin-mixed} is stabilizable, as formalized later in \autoref{sec:main-results}.

Given an initial state $x$, we  are interested in finding sequences of continuous and discrete inputs for system \eqref{eq:sys-lin-mixed} that minimize    cost    $J$. That is, ideally we wish to find for any $x\in\R^{n_x}$ an optimal policy, i.e., the optimal sequences of continuous and discrete inputs giving the \emph{optimal value function}
\begin{equation} \label{eq:Vstar}
    V^\star(x):=\min_{\bm{u},\bm{i}}J(x,\bm{u},\bm{i}), 
\end{equation}
assuming it exists. The optimal policy can be constructed from an optimal feedback law $h^\star$, which is a function obtained from the set-valued map  $H^\star:\R^{n_x}\rightrightarrows \R^{n_u}\times\mathbb{M}$, where
\begin{equation} \label{eq:Hstar}
H^\star: x\mapsto\argmin_{(u,i)}\left\{\ell(x,u,i)+V^\star(A_{i}x+B_{i}u)\right\},
\end{equation}
by making a single choice for the    continuous and discrete input for each possible state, i.e., $h^\star \in H^\star$.

Solving \eqref{eq:Vstar}-\eqref{eq:Hstar}   is known to be  challenging in general, with \cite{weizhang_Value_2009} suggesting that it may be NP-hard due to the combinatorial growth of possible switching sequences. While an exact formulation of the optimal value function remains an open problem, it is known that this function can be approximated arbitrarily closely by using a finite but exponentially growing  set of quadratic functions \cite{zhang_InfiniteHorizon_2012,weizhang_Value_2009}.

To overcome the computational complexity of solving \eqref{eq:Vstar} exactly, we propose a near-optimal point-wise algorithm for the switched LQR 
problem. By focusing on obtaining accurate estimates of the value function \emph{at a given initial state}, rather than 
attempting to solve it across the entire state space, we significantly reduce  the computational burden. In this way, we can 
close the loop of system~\eqref{eq:sys-lin-mixed} with our algorithm, generating inputs that provide desirable guarantees of near-optimality and 
exponential stability for the controlled system. This new algorithm, called \OPlin, is inspired by optimistic planning algorithms in \cite{Hren-Munos-2008optimistic,Busoniu-Munos-12,granzotto_stable_2022}. We provide  explicit relationships between the stability of the system controlled by \OPlin, its near-optimality guarantees, and a tuning parameter related to the maximum amount of computations in a sense given later, under given assumptions. 

The algorithm is presented in the next section after some preliminaries.

\section[\texorpdfstring{\OPlin}{OPlin}]{\OPlin} \label{sec:OPlin}
The algorithm we present approximates $V^\star(x)$ in \eqref{eq:Vstar} at  any given $x\in\R^{n_x}$ by solving finite-horizon optimal control problems with a sufficiently large horizon, with or without a terminal cost.  As we will see in the sequel, the finite-horizon optimal value function is calculated exploiting  time-varying Riccati equations given  sequences of discrete inputs, without the need to calculate the continuous input. This allows for the application of optimistic planning methods \cite{Hren-Munos-2008optimistic,Busoniu-Munos-12,granzotto_stable_2022},  which explore  countable and finite combination of discrete inputs in an efficient manner.

\subsection{Finite-horizon optimal control}

Given $x\in\R^{n_x}$ and a finite horizon $d$, consider the $d$-horizon cost  $J_d$,  defined as
\begin{equation}\label{eq:Jd}
\begin{gathered}
    J_d(x,\bm{u},\bm{i})    :=    \sum_{k=0}^{d-1} \ell(\phi(k,x,\bm{u}_{\llbracket \cdot, k \rrbracket},\bm{i}_{\llbracket \cdot, k \rrbracket}),u_k,i_k)\\
     {}+J_0(\phi(d,x,\bm{u},\bm{i})),
\end{gathered}
\end{equation}
with   $J_0(x):=x^\top \underline{P}x$ for some designed $\underline{P}\succeq0$ of conformable dimensions. We also introduce  the optimal value function for cost \eqref{eq:Jd} for a \emph{given fixed} sequence of discrete inputs $\bm{i}$ of length $d\in\Zo$ and a given $x\in\R^{n_x}$, that is  
\begin{equation} \label{eq:V-of-d-circ}
     V_d^\circ(x,\bm{i}) := \min_{\bm{u}} J_d(x,\bm{u},\bm{i}).
\end{equation}   

The minimization in \eqref{eq:V-of-d-circ} is intentionally performed only over the sequence of continuous inputs in $J_d(x,\cdot,\bm{i})$, while the sequence of discrete inputs $\bm{i}$ remains fixed. This formulation leads to a finite-horizon time-varying LQR problem for which the optimal cost optimal value function satisfies, for any $x\in\R^{n_x}$ and $\bm{i}$ of length $d\in\Zo$, 
\begin{equation} \label{eq:V-of-d-circ-bellman}
     V_d^\circ(x,\bm{i}) =    \min_{u} \left\{\ell(x,u,i_0) +  V_{d-1}^\circ(A_{i_0}x+B_{i_0}u,\bm{i}_{\llbracket 1,\cdot  \rrbracket})\right\},
\end{equation} with $V_0^\circ=J_0$. Importantly, $V_d^\circ$ is given by a quadratic form as formalized in the sequel. To see it,  \textcolor{black}{we introduce the Riccati operator associated with mode $i\in\mathbb{M}$ as $\mathcal{R}_i : \R^{n_x\times n_x}\to \R^{n_x\times n_x}$, which is given by} 
\begin{align*} \label{eq:ricatti-op}
        \MoveEqLeft\mathcal{R}_i(P):={} \numberthis\\ 
        &Q_i + A_i^\top P A_i - A_i^\top P B_i (R_i+B_i^\top P  B_i)^{-1} B_i^\top P A_i
\end{align*}  
  for any $P\in\R^{n_x\times n_x}$. Note that $\mathcal{R}_i(P)\succ0$ whenever $P\succeq0$ as $Q_i,R_i\succ0$ for all $i\in\mathbb{M}$ by a Schur complement test\footnote{Let $\mathcal{P}_i{:=}\begin{bsmallmatrix*}
      A_i & B_i
  \end{bsmallmatrix*}^\top P  \begin{bsmallmatrix*}
      A_i & B_i
  \end{bsmallmatrix*}$ and   $\mathcal{S}_i{:=}\begin{bsmallmatrix*}
      A_i^\top P A_i  +Q_i & A_i^\top P B_i\\
      B_i^\top P A_i  & B_i^\top P B_i +R_i
  \end{bsmallmatrix*}= \mathcal{P}_i+\text{diag}(Q_i,R_i)$, so  the Schur complement of $\mathcal{S}_i$   is $\mathcal{R}_i(P)$. Given $P\succeq0$ and $Q_i,R_i\succ0$, then $\mathcal{P}_i{\succeq}0$, hence $\mathcal{S}_i{\succ}0$  and we conclude that $\mathcal{R}_i(P)\succ0$.}.   The next proposition shows that $V_d^\circ$ in \eqref{eq:V-of-d-circ} can be reduced to a quadratic form involving \eqref{eq:ricatti-op}.
  \begin{prop} \label{prop:Vc-of-i}
Given $d\in\Zp$, for any $d$--length sequence of discrete inputs $\bm{i}$, 
$V^\circ_d(x,\bm{i})= x^\top  P_{\bm{i}} x$
for any $x\in\R^{n_x}$
where  $P_{\bm{i}}\succ0$ is  defined as
\begin{equation} \label{eq:ricatti}
P_{\bm{i}}:=\mathcal{R}_{i_0}\circ\mathcal{R}_{i_1}\circ\cdots\circ\mathcal{R}_{i_{d-1}}(\underline{P}).\end{equation}
Moreover, 
\begin{equation}\label{eq:time-varying-bellman}
V^\circ_d(x,\bm{i})= \ell(x,-K_{\bm{i}}x,{i_0}) +  V^\circ_{d-1}((A_{i_0}-B_{i_0}K_{\bm{i}})x,\bm{i}_{\llbracket 1,\cdot  \rrbracket})
\end{equation}
with $K_{\bm{i}}:=(R_{i_0}+B_{i_0}^\top P_{\bm{i}_{\llbracket 1,\cdot  \rrbracket}} B_{i_0})^{-1}   B_{i_0}^\top   P_{\bm{i}_{\llbracket 1,\cdot  \rrbracket}}   A_{i_0}$.
\mbox{}\hfill$\Box$
  \end{prop} 

\noindent\textbf{Proof:} Let $\bm{i}$ be a sequence of discrete inputs of length $d\in\Zp$.  For all  $x\in\R ^{n_x}$, it follows that   $V_1^\circ(x,\bm{i}_{\llbracket d-1,\cdot\rrbracket})= x^\top \mathcal{R}_{i_{d-1}}(\underline{P})x=\ell(x,-K_{\bm{i}_{\llbracket d-1,\cdot\rrbracket}}x,{i_{d-1}}) +  V^\circ_{0}((A_{i_{d-1}}-B_{i_{d-1}}K_{\bm{i}_{\llbracket d-1,\cdot\rrbracket}})x)$ with $ \mathcal{R}_{i_{d-1}}(\underline{P})\succ0$  and $K_{\bm{i}_{\llbracket d-1,\cdot\rrbracket}}:=(R_{i_{d-1}}+B_{i_{d-1}}^\top \underline{P} B_{i_{d-1}})^{-1}   B_{i_{d-1}}^\top   \underline{P}   A_{i_{d-1}}$, which is obtained  by a Schur complement argument  noting that $Q_i,R_i\,{\succ}\,0$ for all $i\,{\in}\,\mathbb{M}$ and $\underline{P}\,{\succeq}\,0$, see, e.g., \cite[\textcolor{black}{App.} A5.5]{boyd_Convex_2004}. Moreover and for the same reasons, $\ell(x,-K_{\bm{i}_{\llbracket d-1,\cdot\rrbracket}}x,{i_0}) +  V^\circ_{0}((A_{i_{d-1}}-B_{i_{d-1}}K_{\bm{i}_{\llbracket d-1,\cdot\rrbracket}})x)\leq \ell(x,u,{i_{d-1}}) +  V^\circ_{0}(A_{i_{d-1}}x+B_{i_{d-1}}u)$  for all  $x\in\R ^{n_x}$. In view of \eqref{eq:V-of-d-circ-bellman}, we  iterate the above $d-1$ times to obtain \eqref{eq:ricatti} and the proof is complete.  \mbox{}\hfill$\blacksquare$

\autoref{prop:Vc-of-i} shows that $V_{d}^{\circ}$ is given by a quadratic form and that the   first element of the optimal $d$-long sequence of continuous inputs is $-K_{\bm{i}} x$.  In this way, we are able to construct matrices $P_{\bm{i}}\succ0$ and sequences  $\bm{u}$ such that $V_{d}^{\circ}(x,\bm{i})=x^\top P_{\bm{i}} x = J(x,\bm{u},\bm{i})$  given  discrete sequences $\bm{i}$ in a simple and straightforward manner in   \eqref{eq:ricatti}  and \eqref{eq:time-varying-bellman} in terms of \eqref{eq:ricatti-op}.  

So far in this section, the sequence of discrete inputs $\bm{i}$ is fixed, we now relax this requirement and analyze  $V_d^\circ(x,\cdot)$ for any $x\in\R^{n_x}$.

\subsection{Finite-horizon value functions}

\OPlin aims to produce    a monotonically increasing sequence  of estimates of the value function that converge to the optimal value function $V^\star$ in \eqref{eq:Vstar} from below. For this purpose, we make the next assumption on matrix $\underline{P}$, i.e., on the terminal cost in \eqref{eq:Jd}.

\newtheorem{tc}{Condition}
\begin{tc}[\autoref{TC:terminal-cost-general}]\label{TC:terminal-cost-general}   Matrix $\underline{P}\succeq0$ is such that, for any  $i\in\mathbb{M}$,
     \begin{equation} \label{eq:terminal-LMI}
        \begin{bmatrix*}
            A_i^\top \underline{P} A_i -\underline{P}+Q_i & A_i^\top \underline{P} B_i \\ B_i^\top \underline{P} A_i &  R_i + B_i^\top \underline{P} B_i
        \end{bmatrix*}\succeq0.
    \end{equation}
    \mbox{}\hfill$\Box$
\end{tc}

\autoref{TC:terminal-cost-general} is made  without loss of generality as   matrix $\underline{P}$ is a design parameter which we are free to choose, and   there always exists some matrix  $\underline{P}\succeq0$ verifying \eqref{eq:terminal-LMI}. Indeed,    it suffices to take $\underline{P}=0$ as $Q_i,R_i\succ0$ for all $i\in\mathbb{M}$. Other candidates for $\underline{P}$ can be calculated by solving the LMI in \eqref{eq:terminal-LMI}.  \textcolor{black}{For reasons that will become clear later, it is desirable to take  $\underline{P}$ subject to  \eqref{eq:terminal-LMI}  that are ``big'' in some sense, e.g., maximizes $\Tr(\underline{P})$ or maximizes $\text{logdet}(\underline{P})$. subject to  \eqref{eq:terminal-LMI}  that maximizes $\text{tr}(\underline{P})$.} \textcolor{black}{Interestingly, when $M=1$, maximizing $\text{tr}(\underline{P})$ under \eqref{eq:terminal-LMI} leads to the optimal value function for the infinite-horizon standard (non-switched) LQR problem whenever $(A_1,B_1)$  stabilizable and $(A_1,Q_1)$  detectable  \cite{vandenberghe_Semidefinite_1999}.}

\autoref{TC:terminal-cost-general} in conjunction to \autoref{prop:Vc-of-i} leads to the next two statements.

\newcommand\mdoubleplus{\ensuremath{\mathbin{+\mkern-10mu+}}}
\begin{prop} \label{prop:Vc-monotone}
Let  $d\in\Zo$ and $x\in\R^{n_x}$, for any  discrete inputs sequence $\bm{i}$ of length $d$ and any $i\in\mathbb{M}$,
\begin{equation} \label{eq:Vc-monotone}
    V_d^\circ(x,\bm{i})\leq V_{d+1}^\circ(x, \bm{i}\oplus i)
\end{equation}
where $(\bm{i}\oplus  i)_{\llbracket\cdot,d-1\rrbracket}:=\bm{i}$ and $(\bm{i}\oplus  i)_{\llbracket d,d\rrbracket} := i$.
 \mbox{}\hfill$\Box$
\end{prop}

\noindent\textbf{Proof:} Let $x\in\R^{n_x}$, $d\in\Zo$,  sequence of discrete inputs $\bm{i}$ of length $d$ and $i\in\mathbb{M}$.  Furthermore, let $\bm{u}^\star$ and $\bm{\widetilde{u}}^\star$ be   continuous input sequences such that $V_d^\circ(x,\bm{i})=J_d(x,\bm{u}^\star,\bm{i})$ and $V^\circ_{d+1}(x,\bm{i}\oplus i)=J_{d+1}(x,\bm{\widetilde{u}}^\star,\bm{i}\oplus i)$, respectively, which exist by \autoref{prop:Vc-of-i}; note that both $\bm{u}^\star$ and $\bm{\widetilde{u}}^\star$ depend on $\bm{i}$ and $\bm{i}\oplus i$, respectively. In view of \eqref{eq:Jd},  
\begin{equation}
\begin{split}
   \MoveEqLeft[4]  J_{d+1}(x,\bm{\widetilde{u}}^\star,\bm{i}\oplus i)  = J_{d}(x,\bm{\widetilde{u}}^\star_{\llbracket\cdot,d-1\rrbracket},\bm{i})\\   &{}+\ell(\phi_{d},\widetilde{u}^\star_{d},i)+J_0(\phi_{d+1})-J_{0}(\phi_{d}),
\end{split} 
\end{equation}
     where $\phi_{d}:=\phi(d,x,\bm{\widetilde{u}}^\star_{\llbracket \cdot, d-1 \rrbracket},\bm{i})$ and $\phi_{d+1}:= A_{i}\phi_{d}+B_{i}\widetilde{u}^\star_{d}$. Moreover, $(x,u,i)\mapsto\ell(x,u,i)+J_0(A_{i}x+B_{i}u)-J_{0}(x)$ is non-negative  by \autoref{TC:terminal-cost-general}, by pre- and post-multiplying  \eqref{eq:terminal-LMI} with $(x,u)^\top$ and $(x,u)$, respectively.
    Therefore, 
    $V^\circ_{d+1}(x,\bm{i}\oplus i)=J_{d+1}(x,\bm{\widetilde{u}^\star},\bm{i}\oplus i)\geq J_{d}(x,\bm{\widetilde{u}}^{\star}_{\llbracket  \cdot,d-1 \rrbracket},\bm{i})$.   By optimality of $\bm{u}$ vis-a-vis $J(x,\cdot,\bm{i})$, we conclude $V^\circ_{d+1}(x,\bm{i}\oplus i)\geq V^\circ_d(x,\bm{i})$ and the desired result is obtained. \mbox{}$\hfill\blacksquare$

    \autoref{prop:Vc-monotone} shows that   cost  $V^\circ_d$ may only increase in $d$.  As extending a sequence of discrete inputs, i.e., increasing the horizon in \eqref{eq:Jd}, results in a higher or equal cost, we can thus   discard   sequences with high cost and focus on the more promising ones. This fact   will be used in the proposed algorithm as,  when   comparing sequences of different lengths and costs, the sequence with a smaller associated cost $V^\circ$ is preferable to be extended as the larger ones are guaranteed to always produce larger costs as measured by $V^\circ$ by \autoref{prop:Vc-monotone}.   

    The proposed algorithm also relies on the next property of finite-horizon costs \eqref{eq:Jd}.

\begin{prop} \label{prop:Vs-of-d}
    For any $d\in\Zo$ and $x\in\R^{n_x}$, there exists a $d$-length sequence $\bm{i}^{d,\star}$ such that, for any $\bar{d}$-length sequence  $\bm{i}$ with $\bar{d}\geq d$,
\begin{equation}\label{eq:Vd-increase}
 V^\star_d(x):=V^\circ_d(x,\bm{i}^{d,\star})\leq V^\circ_{\overline{d}}(x,\bm{i}).   
\end{equation}
In addition,  
\begin{equation}\label{eq:Vd-lowers-Vstar}
    V^\star_d(x)\leq J(x,\bm{u},\bm{i})
\end{equation}
for any infinite-length sequence of continuous and discrete inputs $\bm{u}$, $\bm{i}$ such that $\lim_{k\to\infty}\phi(k,x,\bm{u},\bm{i})=0$. \mbox{}\hfill$\Box$
\end{prop}
\noindent\textbf{Proof:} Let $x\in\R^{n_x}$ and $\bar{d},d\in\Zo$ such that $\bar{d}\geq d$. For the  existence of  $V^\star_d(x)$ such that $V^\star_d(x)\leq V_d^\circ(x,\bm{i})$ for any $d$-length sequence $\bm{i}$, it suffices to consider $\bm{i}^{d,\star}\in \argmin_{\bm{i}} V^\circ_d(x,\bm{i})$, which exists since $V^\circ_d(x,\bm{i})\geq0$ and the set of $d$-length sequences of discrete inputs is   finite. Moreover, \eqref{eq:Vd-increase} holds for any $\bar{d}$-length sequence $\bm{i}$ where $\bar{d}\geq d$ by   virtue  of \autoref{prop:Vc-monotone}, as necessarily $V^\circ_d(x,\bm{i}_{\llbracket\cdot,d-1\rrbracket})\geq V^\star_d(x)$.

  To show \eqref{eq:Vd-lowers-Vstar},  first notice that it is trivially verified when $J(x,\bm{u} ,\bm{i})=\infty$ as $V_d^\star(x)$ is finite as $V_d^\star(x) \leq x^\top P_{\bm{i}} x$ for any $d$-sequence of discrete inputs $\bm{i}$. Hence, let now $\bm{u}$, $\bm{i}$ be infinite-length sequences such that $\lim_{k\to\infty}\phi(k,x,\bm{u},\bm{i})=0$ and $J(x,\bm{u} ,\bm{i})<\infty$.  Necessarily $J(x,\bm{u} ,\bm{i})=J_{\overline{d}}(x,\bm{u}_{\llbracket \cdot,{\overline{d}}-1\rrbracket},\bm{i}_{\llbracket \cdot,{\overline{d}}-1\rrbracket})+J(\phi_{\overline{d}},\bm{u}_{\llbracket {\overline{d}} ,\cdot\rrbracket},\bm{i}_{\llbracket {\overline{d}} ,\cdot \rrbracket})-J_0(\phi_{\overline{d}})$ where $\phi_{\overline{d}}:=\phi({\overline{d}},x,\bm{u}_{\llbracket \cdot,{\overline{d}}-1\rrbracket},\bm{i}_{\llbracket \cdot,{\overline{d}}-1\rrbracket})$ for any ${\overline{d}}\in\Zo$.  Hence, $J_{\overline{d}}(x,\bm{u}_{\llbracket \cdot,{\overline{d}}-1\rrbracket},\bm{i}_{\llbracket \cdot,{\overline{d}}-1\rrbracket})-J_0(\phi_{\overline{d}})\to J(x,\bm{u} ,\bm{i})$  when ${\overline{d}}\to\infty$, as $J_0(x)=x^\top \underline{P} x$ with $\underline{P}\succeq0$ and $\lim_{k\to\infty}\phi(k,x,\bm{u},\bm{i})=0$, and the proof is concluded in view of $V^\star_d(x)\leq J_{\overline{d}}(x,\bm{u}_{\llbracket \cdot,d-1\rrbracket},\bm{i}_{\llbracket \cdot,d-1\rrbracket}) $ for any ${\overline{d}}\geq d$ as consequence of \autoref{prop:Vc-monotone}.    \mbox{}\hfill $\blacksquare$

\autoref{prop:Vs-of-d} establishes  the existence of an optimal sequence of discrete inputs for cost \eqref{eq:Jd} corresponding to the optimal value function $V_d^\star(x)$, which   is non-decreasing in $d$ by \autoref{prop:Vc-monotone}. Moreover, $V_d^\star$ lower-bounds the incurred cost of any stabilizing sequence of inputs. Hence, $V_d^\star$ also lower bounds the optimal value function for the original infinite-horizon cost \eqref{eq:Vstar}, provided it exists and  solutions to system \eqref{eq:sys-lin-mixed} in closed-loop with $H^\star$ are stabilizing, which will be established later.    

We are ready to describe how a planning approach can efficiently calculate  discrete inputs $\bm{i}^{d,\star}$ and cost $V_d^\star(x)$  given $x\in\R^{n_x}$.

\subsection[\texorpdfstring{The \OPlin algorithm}{The OPlin algorithm}]{The \OPlin algorithm}  

For any given state $x\in\R^{n_x}$ and horizon $d\in\Zo$, \OPlin   explores trees  by prioritizing branches that are most promising based on current cost estimates. In an optimistic sense, the ``best''  candidate  for the infinite-horizon optimal cost is the sequence with the smallest evaluated cost \eqref{eq:Jd}. By exploring the most promising sequences, \OPlin reduces the computational complexity of solving the switched LQR problem while maintaining stability and performance guarantees.   Specifically,  the algorithm   avoids exploring all branches that incur  costs larger than $V_d^\star(x)$.  

The method is given in \autoref{algo:OPlin-budget}. At each step, \OPlin selects for exploration the     sequence of discrete inputs that has the lowest cost calculated so far. This sequence is   extended by adding new discrete inputs, and their cost are  evaluated in terms of    positive semi-definite matrices that solve \eqref{eq:V-of-d-circ}.  The process continues until   horizon $d$ is reached. 

The notation of \autoref{algo:OPlin-budget} is as follows.  We denote by $\mathcal{T}$ the exploration tree from initial state $x\in\R^{n_x}$, constructed from discrete input sequences $\bm{i}$ of length $d$,  with respective cost $V_d^\circ(x,\bm{i})$. A  node $N$ payload (or data) is a triplet  in $ \mathbb{M}^d\times\R^{n_x\times n_x}\times\R_{\geq0}$ for some $d\in\Zo$. These are, respectively: the label of the node which is given by the discrete input sequence, denoted by $\bm{i}(N)$; the quadratic matrix that generates   cost $V_d^\circ(\cdot,\bm{i}(N))$, denoted by $P(N)$; and the point-wise cost $V_d^\circ(x,\bm{i}(N))$, denoted by $J(N)$.   A leaf is a node of $\mathcal{T}$ with no children, and the set of all leaves of $\mathcal{T}$ is denoted $L(\mathcal{T})$. At iteration $j\in\Zo$, a leaf $L_j\in L(\mathcal{T})$ is fully expanded by extending the sequence of discrete inputs of $L_j$, denoted by $\bm{i}(L_j)$. That is, for every $i \in\mathbb{M}$, we add a child to $L_j$ labeled by $\bm{s}^i:=\{\bm{i}(L_j),i\}$,  which are new leaves of $\mathcal{T}$;  such leaves have an associated cost matrix given by $P_{\bm{s}^i}$ given by  \eqref{eq:ricatti-op}-\eqref{eq:ricatti} and point-wise cost given by $x^\top P_{\bm{s}^i} x$; after this, $L_j$ is no longer a leaf, but becomes an inner node. The algorithm repeats this process until it finds a leaf $S$ with a sequence of length $d$, and the computational resources utilized for exploring the tree are denoted as a budget $\mathcal{B}\in\Zp$, which corresponds to $\mathcal{B}$ leaf expansions. Leaf $S$  verifies $\bm{i}^{\star,d}=\bm{i}(S)$ and $V^\star_d(x)=J(S)$, as seen in the next proposition.

\begin{algorithm}[t]
\setlength{\belowdisplayskip}{0pt} \setlength{\belowdisplayshortskip}{0pt}
\setlength{\abovedisplayskip}{0pt} \setlength{\abovedisplayshortskip}{0pt}

\renewcommand{\hypcapspace}{2\baselineskip}
\capstart
\renewcommand{\hypcapspace}{0.5\baselineskip}
 \caption{\strut \label{algo:OPlin-budget} \OPlin}
 \begin{algorithmic}[1]
 \renewcommand{\algorithmicrequire}{\textbf{Input:}}
 \renewcommand{\algorithmicensure}{\textbf{Output:}}
 \Require  State $x$, horizon $d$.
 \Ensure   Sequence $\bm{i}^{d,\star}(x)$, cost $V^\star_d(x)$, budget $\mathcal{B}$.
\State \textbf{Initialization:}
 \Statex $j \gets 0, S\gets \emptyset$ 
 \Statex $\mathcal{T}\gets \{\emptyset,\underline{P},x^\top \underline{P} x\}$   \Comment{sequence and data}
 \Repeat \Comment{tree exploration}
   \Statex $j\gets j+1$
   \State  \textbf{Find a leaf} $L$ that minimize  $J(L)=x^\top P_{\bm{i}(L)} x$
   \begin{equation}\begin{gathered}
      L_j \in \argmin_{L\in\mathcal{T}} J(L) \\
      \mathcal{T}\gets\mathcal{T}\setminus L_j
   \end{gathered}\tag{OP.1}
  \end{equation}
  \State    \textbf{Add} the children of $L_j$ to $\mathcal{T}$
  \Statex \qquad   For all $i\in\mathbb{M}$ and $\bm{s}^i :=\{\bm{i}(L_j),i\}$,  
\begin{equation}
\begin{gathered}    
P_{\bm{s}^i}:=
\mathcal{R}_{s^i_0}\left (P_{\bm{s}^i_{\llbracket  1,\cdot  \rrbracket}}\right )\\
  \mathcal{T}\gets\mathcal{T} \cup \{ \bm{s}^i,\ P_{\bm{s}^i},\ x^\top P_{\bm{s}^i} x \}
 \end{gathered} \tag{OP.2} \label{OP:ith-leaf} 
 \end{equation}
\hypertarget{OP:Select}{\Until{$\text{length}(\bm{i}(L_j))=d$}  \Comment{stopping criterion}}
 \State $S\gets L_j$, \ \Return   $\bm{i}(S)$, $J(S)$, $\mathcal{B}\gets j. ~\mbox{}\hfill$ 
 \end{algorithmic}

\end{algorithm}

\begin{prop} \label{prop:OPlin-basic}
    For every $d\in\Zo$ and $x\in\R^{n_x}$,  \autoref{algo:OPlin-budget} terminates with finite budget $\mathcal{B}\leq \mathcal{B}^\star:=\frac{M^{d+1}-1}{M-1}+1$.  Moreover, $i^{d,\star}(x)=\bm{i}(S)$ and $V^\star_d(x) = J(S)$ hold.\mbox{}\hfill $\Box$
\end{prop}
\noindent\textbf{Proof:}  The proof is adapted  from \cite[Proof of Proposition 2]{granzotto_stable_2022} to cope with terminal cost $J_0\neq0$ and a given horizon $d$. The termination of \autoref{algo:OPlin-budget}   follows from the fact that any tree with a leaf at  depth $d$ has at most the same number of nodes than the shallowest (and densest) tree of depth $d$. In turn, the shallowest tree  with depth $d$ requires $1+M+\cdots+M^{d}+1=\frac{1-M^{d+1}}{1-M}+1=\mathcal{B}^\star$ leaf expansions.  Therefore, a leaf with a sequence of length $d$ must have been selected for some $j \leq \mathcal{B}^\star$, hence the \hyperlink{OP:Select}{stopping criterion}   is eventually verified for $j=\mathcal{B}\leq \mathcal{B}^\star$.
Similarly,  $J(S)=V_d^\star(x)$ is implied by the fact that the first sequence selected by \OPlin with length $d$ is necessarily smaller or equal than any other sequence of length $d$. Indeed, if this was not the case, there would exist a node $N$, not in tree $\mathcal{T}$, such that $J(N)<J(S)$ and $\bm{i}(N)$ of length $d$. Moreover, since $N\not\in\mathcal{T}$, there exists some  leaf $L\in\mathcal{T}$ whose sequence $\bm{i}(L)$ extends to $\bm{i}(N)$; or, in other words, where $N$ is a descendant of $L$. In view of \autoref{prop:Vc-monotone}, necessarily $J(L)\leq J(N)$. But $J(S)\leq J(L)$ for any leaf $L\in\mathcal{T}$ at step $j=\mathcal{B}$ in view of \eqref{OP:ith-leaf} and  the \hyperlink{OP:Select}{stopping criterion},  then $J(N)<J(S)\leq J(L)\leq J(N)$ and a contradiction is attained.  Thus $J(S)\leq J(N)$ for any $N$, and since  $\bm{i}(N)$ is an arbitrary sequence of length $d$, $J(S) \leq V_d^\star(x)$. Equality is verified as $J(S) \geq V_d^\star(x)$ holds in view of $\bm{i}$ being a $d$-length sequence and \autoref{prop:Vs-of-d}.   \mbox{}$\hfill\blacksquare$

\autoref{prop:OPlin-basic} implies that \OPlin solves  the optimal finite-horizon problem \eqref{eq:Vd-increase}, and terminates within at most $\mathcal{B}^\star$ iterations. Note that  $\mathcal{B}^\star$ is a very conservative estimate based on the worst case exploration, in which case \OPlin behaves as a brute-force algorithm. In practice, \OPlin typically requires  much fewer iterations by only exploring  a few branches,   see \autoref{sec:example} for examples.

In the remainder of the paper, we relate the finite-horizon optimal control problem solved by \OPlin to the original infinite-horizon optimal control problem solved by \autoref{sec:notation}, in a way that allows us to determine a suitable $d$ with explicit near-optimality and stability guarantees. The corresponding analysis is inspired by \cite{granzotto2020finite,granzotto_stable_2022} and extends to the case of a non-zero terminal cost. 

\begin{rem}
For clarity,  \OPlin \!\! iterates   until  the $d$-horizon optimal sequence is found, allowing the budget $\mathcal{B}$ to adapt. Alternatively, fixing $\mathcal{B}$  yields similar guarantees, provided $\mathcal{B}$ is large enough to find long enough sequences.
\mbox{}\hfill$\Box$  
\end{rem}

\begin{rem}
    There are two   implementation details that can be  employed to reduce computational complexity of \autoref{algo:OPlin-budget}. The first  is to  employ  a sorted list  (or a \emph{priority queue}) of $J(L)$ for all $L\in\mathcal{T}$, which allows to quickly obtain  the leaf $L\in\mathcal{T}$ that minimizes $J(L)$. This is relevant when the size of the tree grows large, avoiding a search over all leaves to find the one with smallest cost at every iteration of \OPlin.  The second  is that $P_{\bm{i}}$ is recursively defined and hence may be \emph{memoized} to avoid recomputations when two or more  sequences share   subsquences. Indeed, to calculate $P_{\bm{i}}$ for some sequence $\bm{i}$ of length $d$, we also need to calculate $P_{\bm{i}_{\llbracket  k, \rrbracket}}$ for all $k\in\{1,\ldots,d-1\}$, which may be already been calculated.  By employing  a cache of all previously computed $P_{\bm{i}}$ we  avoid  repeated calculations of $\mathcal{R}_{i_0}(P_{\bm{s}^i_{\llbracket  1,\cdot \rrbracket}})$, which is relevant when $n_x$ is very large and \eqref{eq:ricatti} in \eqref{OP:ith-leaf} is an expensive calculation. For similar reasons, we can reuse such cache of matrices to any $x\in\R^{n_x}$.  This is a  key advantage  of \OPlin compared to  \cite{Hren-Munos-2008optimistic,Busoniu-Munos-12,granzotto_stable_2022}, as we do not look into future states to obtain the incurred cost of a sequence, instead we calculate the positive
semi-definite matrices that solve \eqref{eq:V-of-d-circ}, which  ascertains the
cost of the sequence independently of $x$. Note though that while these matrices produces a cost for any $x$, \OPlin must be run for each $x$ to find the correct    discrete sequence and associated cost.  In simulations later, we will employ  the priority queue but not the cache.~\mbox{}\hfill$\Box$
\end{rem}

\section{Guarantees}  \label{sec:main-results}

In the previous section, we established that \OPlin solves the optimal finite-horizon problem for cost \eqref{eq:Jd}. We now establish  \OPlin as a tool to approximate $V^\star$. In addition to the requirement on $\underline{P}$ in \autoref{TC:terminal-cost-general}, we require for this purpose some extra but mild assumptions.

\subsection{Stabilizability assumption}
The next assumption  may be verified by constructing inputs that drive the system to the origin sufficiently fast and with quadratic upper-bound on the incurred cost.
\begin{sass}[\autoref{SA:well-posed}]\label{SA:well-posed} 
There exist a real matrix  $\overline{P}\succ0$ and  some sequence of inputs $\bm{u}(x)$ and $\bm{i}(x)$ such that   $V^\star(x)=J(x,\bm{u}(x),\bm{i}(x)) \leq x^\top \overline{P} x$  for any $x\in\R^{n_x}$. \mbox{}\hfill $\Box$
\end{sass}

\textcolor{black}{Due to space constraints, we omit the verification of the conditions in \cite{Keerthi-Gilbert-tac85} that guarantee the existence of optimal inputs.} Given the existence of $V^\star$, a sufficient condition for the upper-bound  is  $(A_i,B_i)$  stabilizable for at least one $i\in\mathbb{M}$. In this case, a suitable  $\overline{P}$  is obtained by solving $A_\text{cl}^\top \overline{P} A_\text{cl}{-} \overline{P}{=}{-}Q_i{-}(K^0)^\top R_i K^0$ with $A_\text{cl}:=(A_i{-}B_iK^0)$  for $K^0\in\R^{n_u\times n_x}$ such that $A_i-B_iK^0$ is Schur. 

\autoref{SA:well-posed}    implies  the next property on the optimal value function and optimal inputs. The proof is omitted as it is a direct application of Propositions \ref{prop:Vs-of-d} and \ref{prop:OPlin-basic}.
\begin{lem} \label{lem:Vstar-bound}
    For any $x\in\R^{n_x}$, 
    \begin{equation} \label{eq:Vd-Vs-lower-bound}
        V^\star_d(x)\leq V^\star(x) \leq x^\top \overline{P} x.
    \end{equation}
    Moreover, $\lim_{k\to\infty}\phi(k,x,\bm{u}^\star(x),\bm{i}^\star(x))=0$  for sequences $\bm{u}^\star(x)$ and $\bm{i}^\star(x)$ such that $V^\star(x)=J(x,\bm{u}^\star(x),\bm{i}^\star(x))$.\mbox{}\hfill$\Box$
\end{lem} 

\autoref{lem:Vstar-bound} ensures that $V^\star_d$ calculated by \autoref{algo:OPlin-budget} lower bounds  the optimal value function $V^\star$, and moreover that solutions to \eqref{eq:sys-lin-mixed} in closed-loop with optimal inputs converge to the origin. 

We now define key constants that characterizes the stability  and near-optimality properties provided by  \OPlin.

\subsection{Key constants}

\textcolor{black}{In the remainder of this section, we unravel the relationship of both the horizon~$d$ and the initial cost function $\underline{P}$ to the stabilizing and the induced near-optimality properties of \OPlin.}  Specifically, we show that large  horizons lead to exponential convergence  to the origin of any solution to the corresponding closed-loop system, as well as  tighter near-optimality bounds.  Moreover, the required horizon decreases when the initial value functions are closer to the optimal value. 

To make explicit the above relationship, we introduce $\alpha \in (0,1)$ and $\alpha_{0}>0$  sufficiently small  verifying, for all $i\in\mathbb{M}$,
 \begin{equation} \label{eq:constants-alpha}
 \begin{gathered}
       \alpha\overline{P}\preceq Q_i\qquad \alpha_{0}(\overline{P}-\underline{P})\preceq Q_i.           
 \end{gathered}
 \end{equation}
  \textcolor{black}{As  $\overline{P},Q_i\succ0$ and $\underline{P}\succeq0$ are known matrices, the only unknowns  in \eqref{eq:constants-alpha} are the scalars
$\alpha$ and $\alpha_0$. 
Such inequalities are always feasible, since the choice
$\alpha=\alpha_{0}=\nicefrac{\min\limits_{i}{\underline{\lambda}(Q_i)}}{\overline{\lambda}(\overline{P})}$ satisfies \eqref{eq:constants-alpha}. This choice provides a simple, explicit, and conservative lower bound, while the maximum feasible values of $\alpha$ and $\alpha_{0}$ can be efficiently computed using LMI solvers.} It is important to note that $\alpha$ is related to the guaranteed decay rate\footnote{This will become clear later, see \eqref{eq:cl-stab} with $d=\infty$.} to the origin of solutions to system \eqref{eq:sys-lin-mixed} in closed-loop with $H^\star$ in \eqref{eq:Hstar}, while $\alpha_{0}$ is related to the initial near-optimality gap\footnote{\textcolor{black}{Hence why it may be desirable to take $\underline{P}$ ``large''.}}, as $V^\star_0(x)=x^\top \underline{P} x\leq V^\star(x)\leq x^\top  \overline{P} x$ under \autoref{SA:well-posed}. That is, larger $\alpha$ is related to faster exponential decay,  while larger $\alpha_0$ is related to smaller initial near-optimality gap.  Together, they dictate how large $d$ must be  in order to provide stability and near-optimality properties for \OPlin, as we   establish next.

\subsection{Stability guarantees}
We consider system \eqref{eq:sys-lin-mixed} in closed-loop with \autoref{algo:OPlin-budget} with a given horizon $d\in\Zo$ for all $x\in\R^{n_x}$, that is, 
\begin{equation} \label{eq:cl-sys}
    x^+\in \left\{A_i x +B_i u \mid (u,i)\in H_d^\star(x) \right \}=:F_d(x)
\end{equation}
for all $x\in\R^{n_x}$ where $H_d^\star(x):=\bigl\{  (u,i)\in\R^{n_u}\times\mathbb{M} \mid \allowbreak V^\star_{d}(x)  =\ell(x,u,i)+V^\star_{d-1}(A_{i}x+B_{i}u)\bigr\}$. Note that $(-K^\star_{\bm{i}(S)} x,\bm{i}_0(S))\in H_d^\star(x)$ with  $\bm{i}(S)$ generated by \autoref{algo:OPlin-budget} and $K^\star_{\bm{i}(S)}$   defined in \autoref{prop:Vc-of-i}. The notation $\phi^\star_d$ stands for solutions to \eqref{eq:cl-sys}, i.e., $\phi^\star_d$ verify $\phi^\star_d(k+1,x)\in F_d(\phi^\star_d(k,x))$ and $\phi^\star_d(0,x)=x$ for any $x\in\R^{n_x}$ and $k\in\Zo$.

We establish global exponential stability  for  system \eqref{eq:cl-sys} provided $d$ is sufficiently large to guarantee exponential decrease.  
\begin{thm} \label{thm:cl-stab}
 Let $d>\max\left\{1,\frac{\log{\alpha_{0}\alpha}}{\log{1-\alpha}}+1\right\}$ with $\alpha$ and $ \alpha_0$ verifying \eqref{eq:constants-alpha}. Then,  for any $x\in\R^{n_x}$,
 \begin{equation} \label{eq:cl-stab}
     \left | \phi^\star_d(k,x)\right) |  \leq   \beta \lambda_d^k   |x|
 \end{equation}
   with $\lambda_d:=1-\alpha+\frac{(1-\alpha)^{d-1}}{\alpha_{0}} \in(0,1)$ and  $\beta:=\nicefrac{\overline{\lambda}(\overline{P})}{\min_{i}{\underline{\lambda}(Q_i)}}$. \mbox{}\hfill$\Box$
\end{thm}
 \noindent\textbf{Proof:} Let $d>1$ and $x\in\R^{n_x}$. We will show that $V^\star$ is a Lyapunov function for \eqref{eq:cl-sys}. To  lower bound $V^\star(x)$, it suffices to recall $\alpha\overline{P}\preceq Q_i$ for all $i\in\mathbb{M}$ given \eqref{eq:constants-alpha}, hence $x^\top \alpha\overline{P}x  \leq V^\star(x) \leq x^\top \overline{P} x$,  where we invoke \autoref{lem:Vstar-bound} for an upper-bound. We now show that system \eqref{eq:cl-sys} verifies an input-to-state like property with respect to  function $V^\star$ and horizon $d$. For this purpose, we build an infinite sequence of inputs such that the first $d-1$ elements corresponds to  the inputs that attain cost $V^\star_d(x)$. That is, let infinite length $\bm{i}$ and $\bm{u}$ be       such that $\upsilon=A_{i_0}x+B_{i_0}u_0$ for some $\upsilon\in F_d(x)$ and $\bm{u}_{\llbracket \cdot, d-2\rrbracket}=\bm{u}^{d,\star}_{\llbracket \cdot, d-2\rrbracket}, \bm{i}_{\llbracket \cdot, d-2\rrbracket}=\bm{i}^{d,\star}_{\llbracket \cdot, d-2\rrbracket}$  with $V^\star_d(x)=J_d(x,\bm{u}^{d,\star},\bm{i}^{d,\star})$. Note that $(u_0,i_0)\in H^\star_d(x)$. We denote the solutions to \eqref{eq:sys-lin-mixed} to inputs $\bm{u}$ and $\bm{i}$ by $x_k$ for any $k\in\Zo$, i.e., $x_k:=A_{i_{k-1}}x_{k-1}+B_{i_{k-1}}u_{k-1}$ with $x_0:=x$ and denote the solutions to \eqref{eq:sys-lin-mixed} to inputs $\bm{u}^{d,\star}$ and $\bm{i}^{d,\star}$ by $\phi_k$ for any $k\in\{0,\ldots,d\}$, note thus:  $x_k=\phi_k$ for $k\in\{0,\ldots,d-1\}$ by construction;  $x_{d}$ may or may not be equal $\phi_d$; only $x_k$ is defined for $k>d$. Thus $x_1=\phi_1=\upsilon$. Without loss of generality, we assume that the inputs for time step $k\geq d-1$ are optimal for the infinite-horizon cost \eqref{eq:J} at state $x_{d-1}$, i.e., $V^\star(x_{d-1})=J(x_{d-1},\bm{u}_{\llbracket d-1, \cdot \rrbracket},\bm{i}_{\llbracket d-1,\cdot\rrbracket})$. By definitions of $V^\star(\upsilon)$, $V_{d-1}^\star(x_1)$, costs $J$ and $J_d$ in \eqref{eq:J} and \eqref{eq:Jd}, we have 
\begin{align*}
        V^\star(\upsilon) &{}\leq   J(x_1,\bm{u}_{\llbracket 1, \cdot \rrbracket}
        ,\bm{i}_{\llbracket 1, \cdot \rrbracket})  \\
        &{}= \textstyle{\sum_{k=1}^{d-2}} \ell(x_k,u_k,i_k)+\sum_{k=d-1}^{\infty} \ell(x_k,u_k,i_k)\\
        &{}= \textstyle{\sum_{k=1}^{d-2}} \ell(x_k,u_k,i_k)+V^\star(x_{d-1})\\
        &\qquad{}\pm\ell(x_{d-1},u^{d,\star}_{d-1},i^{d,\star}_{d-1})
        \pm J_0(\phi_d)\\
        &{}= V^\star_{d-1}(x_1)\numberthis\label{eq:Vs-upsilon-breaking-sequence}\\&\quad{}-\ell(x_{d-1},u^{d,\star}_{d-1},i^{d,\star}_{d-1})-J_0(\phi_{d})+V^\star(x_{d-1}).   
    \end{align*}
In view of \autoref{TC:terminal-cost-general}, 
$0\leq\ell(\phi_{d-1},u^{d,\star}_{d-1},i^{d,\star}_{d-1})+J_0(\phi_d)-J_0(\phi_{d-1})$ as $\phi_d=A_{i^{d,\star}_{d-1}} \phi_{d-1}+B_{i^{d,\star}_{d-1}} u^{d,\star}_{d-1}$, and since $\phi_{d-1}=x_{d-1}$, 
we have $-\ell(\phi_{d-1},u^{d,\star}_{d-1},i^{d,\star}_{d-1})-J_0(\phi_d)\leq -J_0(x_{d-1})$.
Hence, combining the latter with \eqref{eq:Vs-upsilon-breaking-sequence},
\begin{equation}
V^\star(\upsilon) \leq    
         V^\star_{d-1}(x_1)+(V^\star-J_0)(x_{d-1}),
\end{equation}
thus, given that $(u_0,i_0)\in H^{\star}_{d}(x)$,
\begin{equation}\label{eq:cost-decrease}
    V^\star(\upsilon) \leq V^\star_d(x)-\ell(x,u_0,i_0) +(V^\star-J_0)(x_{d-1}),
\end{equation} 
which, by \eqref{eq:Vd-Vs-lower-bound}, becomes
\begin{equation}\label{eq:lyap-temp}
\begin{split}
    V^\star(\upsilon)&{}\leq V^\star_d(x)-x^\top\alpha\overline{P}x +(V^\star-J_0)(x_{d-1})\\
    &{}\leq V^\star(x)-\alpha V^\star(x) +(V^\star-J_0)(x_{d-1}).
\end{split}
\end{equation} 
We now bound $(V^\star-J_0)(x_{d-1})$. By construction of $x_k$ for $k\in\{0,\ldots,d-1\}$ and \eqref{eq:time-varying-bellman} in \autoref{prop:Vc-of-i},  
\begin{equation}
    V^\star_{d-k}(x_k) = \ell(x_k,u_k,i_k) + V^\star_{d-k-1}(x_{k+1}),
\end{equation}
hence,  for all $x\in\R^{n_x}$ and $k\in\{0,\ldots,d-2\}$,
\begin{equation} \label{eq:Y-open-decrease-temp}
     V^\star_{d-k-1}(x_{k+1}) \leq -x_k^\top \alpha\overline{P} x_k +  V^\star_{d-k}(x_k).
\end{equation}
 Note that, similar to the bounds of $V^\star(x)$, $x^\top \alpha\overline{P} x\leq V^\star_{j}(x)\leq V^\star(x) \leq x^\top\overline{P}x$ as $V^\star_j\leq V^\star$ for any $j\in\Zo$. Therefore,  $\alpha  V^\star_{j}(x) \leq x^\top \alpha\overline{P} x$  holds for all $x\in\R^{n_x}$ with $\alpha$ verifying \eqref{eq:constants-alpha}. In particular, $-x_k^\top \alpha\overline{P} x_k \leq -\alpha V^\star_{d-k}(x_k)$, hence in view of \eqref{eq:Y-open-decrease-temp} and for all $k\in\{0,\ldots,d-2\}$,
\begin{equation} \label{eq:Y-open-decrease}
    V^\star_{d-k-1}(x_{k+1}) \leq (1-\alpha) V^\star_{d-k}(x_k).
\end{equation}
By iterating  \eqref{eq:Y-open-decrease} for $d-1$ times and as $V^\star_{1}(x_{d-1})\geq  x_{d-1}^\top Q_{i_{d-1}} x_{d-1}$, we obtain $(1-\alpha)^{d-1} V^\star_{d}(x) \geq x_{d-1}^\top Q_i x_{d-1}$ for any $i\in\mathbb{M}$. Considering \eqref{eq:Vd-Vs-lower-bound} and $\alpha_{0}$ as in \eqref{eq:constants-alpha},  we derive $(1-\alpha)^{d-1} V^\star_{d}(x)\geq \alpha_{0} (V^\star-J_0)(x_{d-1})$, hence
\begin{equation}\label{eq:Y-last-error-bound}
(V^\star-J_0)(x_{d-1})\leq \nicefrac{(1-\alpha)^{d-1}}{\alpha_{0}} V^\star(x).    
\end{equation} By employing the latter  into \eqref{eq:lyap-temp}, we have
\begin{equation}
    V^\star(\upsilon)\leq (1-\alpha+\nicefrac{(1-\alpha)^{d-1}}{\alpha_{0}}) V^\star(x).
\end{equation} 
By iterating the above and  since $\min_{i}\{\underline{\lambda}(Q_i)\} |x|^2  \leq V^\star(x) \leq \overline{\lambda}(\overline{P})  |x|^2$ for all $x\in\R^{n_x}$,
\begin{equation} 
     |\phi^\star_d(k,x)|^2\leq \frac{ \overline{\lambda}(\overline{P})}{\min_{i} \underline{\lambda}(Q_i)} (1-\alpha+\nicefrac{(1-\alpha)^{d-1}}{\alpha_{0}})^k V^\star(x)
\end{equation}
for all $k\in\Zo$. The proof is concluded with $d$ sufficiently large in order to guarantee $(1-\alpha+\nicefrac{(1-\alpha)^{d-1}}{\alpha_{0}})<1$. \mbox{}\hfill$\blacksquare$

\autoref{thm:cl-stab}  shows that global uniform exponential stability of the origin is achieved provided that the horizon $d$ is sufficiently large. Moreover, the decay rate $\lambda_d$ decreases up  to $1-\alpha$  as $d\to\infty$. This theorem thus provides a guaranteed upper-bound on the exponential decay rate of solutions to system~\eqref{eq:cl-sys}, which can be made  as close as desired to the nominal guaranteed decay rate $1-\alpha$, at the cost of more computations.  The required number of computations not only depends on the decay rate estimate $1-\alpha$, but also on the initial mismatch between the optimal value function and cost $V_0^\star$, captured by $\overline{P}-\underline{P}$. Indeed, as $\overline{P},\underline{P}\to V^\star$, $\alpha_{0}\to\infty$. This implies that  good initial   bounds on $V^\star$ reduce the required number of iterations for stability, as desired.

We now provide near-optimal guarantees for the obtained cost $V^\star_d$ and for system \eqref{eq:cl-sys}.

\subsection{Near-optimality guarantees}
The next result upper-bounds the gap between the optimal cost and the cost computed by \OPlin, and such upper-bound is shown to decrease exponentially fast as we increase the horizon $d$. \textcolor{black}{Thus, increasing $d$ leads to increasingly accurate approximations of the optimal value function.}
\begin{thm}\label{thm:near-optimal}
For any $d>1$ and   $x\in\R^{n_x}$, 
\begin{equation} \label{eq:open-loop-near-optimal}
    V^\star(x)-V^\star_d(x)\leq \nicefrac{1}{\alpha_{0}} (1-\alpha)^{d-1}  x^\top \overline{P} x
\end{equation}  holds with $\alpha\in(0,1)$   and $\alpha_{0}$ as in \eqref{eq:constants-alpha}. \mbox{}\hfill$\Box$
\end{thm}
\noindent\textbf{Sketch of proof:} Let $x\in\R^{n_x}$ and $d>1$.  By \eqref{eq:cost-decrease}, we have
 $V^\star(x) \leq V_d^\star(x)+ (V^\star-J_0)(x_{d-1})$.
The proof is concluded by employing the bound  \eqref{eq:Y-last-error-bound}.
\mbox{}\hfill$\blacksquare$

In \autoref{thm:near-optimal}, we see that the cost computed by \autoref{algo:OPlin-budget} approaches $V^\star$ as close as desired, provided we take $d$ sufficiently large. Moreover and similarly to \autoref{thm:cl-stab}, the near-optimality bound is small when the initial
``uncertainty'' captured by $\alpha_0$  also is. Indeed, when $\alpha_{0}\to\infty$, the error bound in \eqref{eq:open-loop-near-optimal}  goes to zero.

\autoref{thm:near-optimal} allows for the next relative near-optimality property, and follows for the same reasons as \autoref{thm:near-optimal}, hence the proof is omitted. 

\begin{prop}
For any $d>1$ and   $x\in\R^{n_x}\setminus\{0\}$, 
\begin{equation} \label{eq:open-loop-near-optimal-relative}
    \frac{V^\star(x)-V^\star_d(x)}{V^\star(x)}\leq \nicefrac{1}{\alpha_{0}} (1-\alpha)^{d-1}  
\end{equation}  holds with $\alpha\in(0,1)$   and $\alpha_{0}$ as in \eqref{eq:constants-alpha}. \mbox{}\hfill$\Box$
\end{prop}

The above results demonstrate that \OPlin  produces inputs with  stability and near-optimality properties, provided that the horizon is large enough for the chosen terminal cost. 
We illustrate the computational aspects of \autoref{algo:OPlin-budget} in an example in the next section.
 

\section{Example} \label{sec:example}

 We consider the example in 
\cite{zhang_InfiniteHorizon_2012} with $n_x=2,n_u=1$ and $M=2$, whose data is given in \autoref{table:data}.  We take $Q_1=Q_2=\mathbb{I}_2$ and $R_1=R_2=1$. Since $(A_1,B_1)$ is controllable, \autoref{SA:well-posed} holds and we construct $\overline{P}$  by solving the LQR problem associated with mode $i=1$.  Thus, we estimate $\alpha \approx 0.14$. All that remains is to find $\underline{P}$ that verifies \autoref{TC:terminal-cost-general}, which we do so by  maximizing $\text{tr}(\underline{P})$ subject to \eqref{eq:terminal-LMI} via an LMI. We find that $\alpha_0 \approx 0.53$ \textcolor{black}{by solving the LMI corresponding to \eqref{eq:constants-alpha}}. Hence, we must employ    $d>18$    to invoke \autoref{thm:cl-stab} stability guarantees, for which a  brute-force approach requires at least $2^{20}-1\approx10^6$ iterations.

\begin{figure}
\centering
\begin{subfigure}{0.5\linewidth}
      \centering
    \includegraphics[width=1\linewidth]{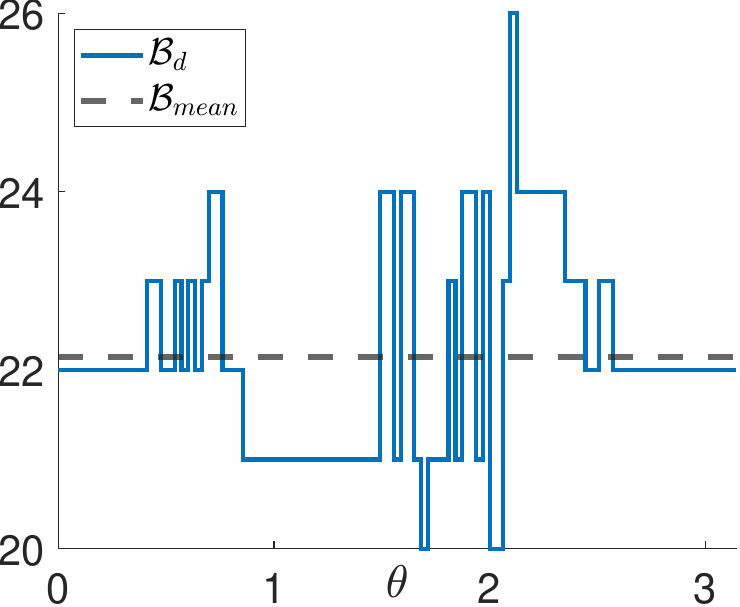}
    \vspace{-4ex}
    \caption{Budgets.}\vspace{-1ex}
    \label{fig:budgets}

\end{subfigure}%
\begin{subfigure}{0.5\linewidth}
    \centering
    \includegraphics[width=1\linewidth,trim={0 0 0 1cm},clip]{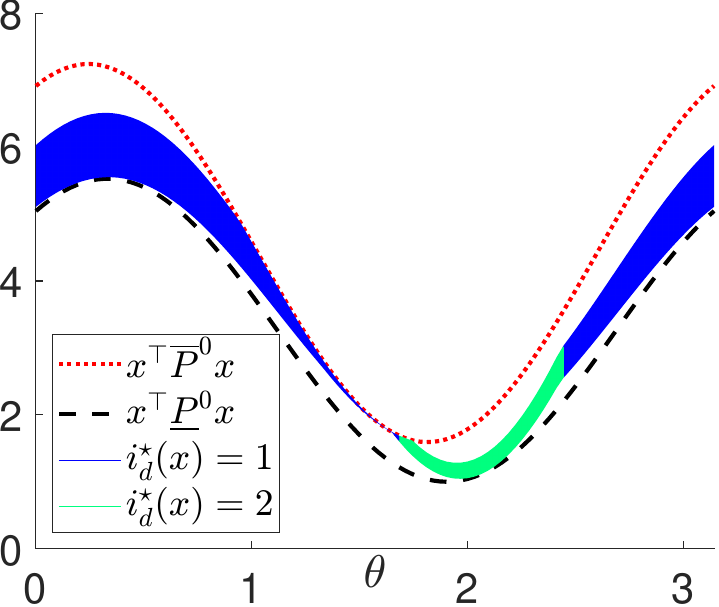}
    \vspace{-4ex}
    \caption{Value function.}\vspace{-1ex}
    \label{fig:values}
\end{subfigure}
\caption{\scriptsize Respectively, Budgets and the Value function for $d{=}19$ and $x\in\{(cos(\theta),sin(\theta)),\theta\in(0,\pi)\}$, where colors blue and green  indicates $i_d^\star=1$ and  2, respectively. $V^\star(x)$ is guaranteed to be inside the blue and green patches.}
\end{figure}

\begin{figure}[t]
    \centering
    \includegraphics[width=0.45\linewidth]{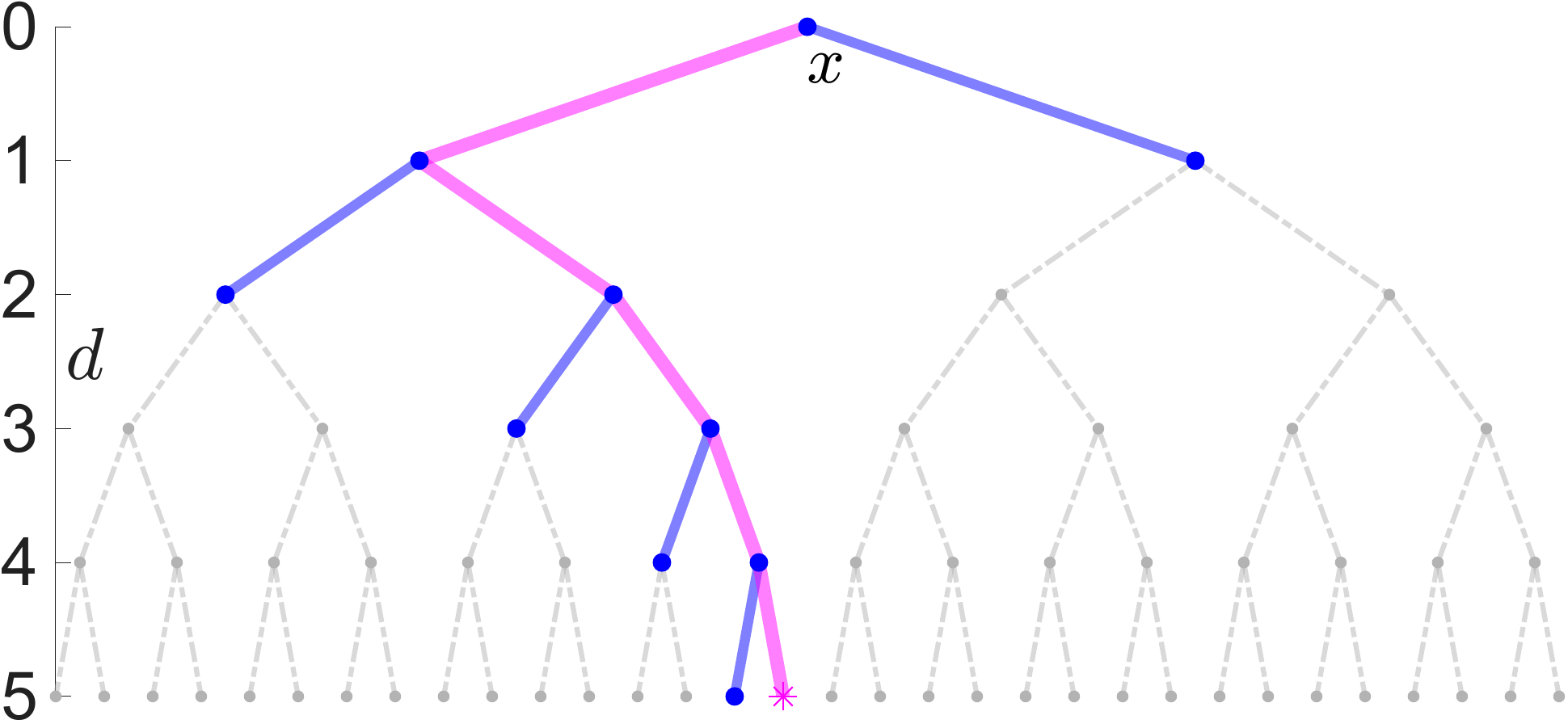}\vspace{-1ex}
    \caption{\scriptsize The brute-force tree for $d=5$. For $x=(-1,0)$, \OPlin explores the tree $\mathcal
{T}$ in blue  and selects leaf $\star$ whose  sequence is in magenta.}\vspace{-0.5ex}
    \label{fig:tree} 
\end{figure}

We apply \OPlin to the example with $d=19$. 
The associated budgets are shown in \autoref{fig:budgets}, and we see $\mathcal{B}_\text{mean}\approx 22$ and $\mathcal{B}_\text{max}= 26$, and a typical  tree developed by \OPlin is found in \autoref{fig:tree}. It appears  that the horizon scales almost linearly with the number of iterations in this example.  This attests the massive  reduction in computations to reach high-depth sequences, as a breadth-first approach would require $10^6$ iterations to reach a similar horizon. We plot in \autoref{fig:values} the initial lower and upper bounds of $V^\star$ given by $\underline{P}$ and $\overline{P}$, as well as the theoretical guaranteed bounds for $V^\star$ given the calculated $V^\star_d$ with $d=19$ by \autoref{thm:near-optimal}. Interestingly, \OPlin actually finds  the optimal value function up to machine precision for all points in this example for horizon $d\leq15$. Indeed, for some $d(x)\leq 15$ we numerically obtain $(V^\star_{d(x)}-V_{d(x)-1}^\star)(x)\approx -10^{-15}$, hence $V^\star_{d(x)}(x)<V_{d(x)-1}^\star(x)$, which is not possible in view of \autoref{prop:Vs-of-d}.  This indicates that the lower bound on $d$ incurs significant conservatism.  We will investigate such conservatism in future work by employing a fine-tuned stopping criterion similar to \cite{granzotto_When_2021}. 


\begin{table}[t]

{\allowdisplaybreaks\scriptsize 
\begin{gather*}
\toprule
\begin{aligned}
    A_1 &{}:= 
    \begin{bmatrix*}[r]
      2 &   1 \\
      0 &  1 \\
    \end{bmatrix*} &   
        A_2&{}:= 
    \begin{bmatrix*}[r]
       2 &   1 \\
      0 &   \nicefrac{1}{2} \\
    \end{bmatrix*} \\   
B_1&{}:= \begin{bmatrix*}[r]  1 &  1\\\end{bmatrix*}^\top   &
B_2&{}:= \begin{bmatrix*}[r]  1 &  2\\\end{bmatrix*}^\top \end{aligned}\\
\begin{aligned}
\overline{P}:=  \begin{bmatrix*}[r]
        6.91 &  1.32\\
       1.32   &  1.92
\end{bmatrix*} &   &
\underline{P}:=  \begin{bmatrix*}[r]
        5.05 &  1.40\\
      1.40   &  1.5
\end{bmatrix*}
\end{aligned}\\\bottomrule 
\end{gather*}}
\vspace{-6ex}
\caption{\label{table:data}}
\end{table}

\vspace{-1ex}
\section{Conclusion} \label{sec:conc}
\vspace{-1ex}

By balancing computational efficiency with near-optimal performance and stability, \OPlin  demonstrates potential for application  in  fields such as robotics, automotive systems, and aerospace, where both stability and near-optimality are essential. A key advantage of \OPlin over ad hoc methods is its ability to compute the value function with known bounds, ensuring near-optimality and providing a point of reference for evaluating heuristic methods lacking formal performance guarantees.

In future work, we plan to reduce  the conservatism in horizon $d$ for stability purposes by designing a suitable stopping criterion, as well as to relax the condition that $Q_i\allowbreak{\succ}0$ for any $i\in\mathbb{M}$ and replacing it with a detectability condition.

\bibliographystyle{IEEEtranS}
\bibliography{IEEEabrv,bib_global,OPlin}

\end{document}